\documentclass [a4,11pt]{amsart}
\usepackage{amsmath,amssymb,amsfonts,enumerate,amsthm, amscd,}

\newcommand{\N}{\mathbb{N}}

\newtheorem{thm}{Theorem}[section]
\newtheorem{cor}[thm]{Corollary}
\newtheorem{lem}[thm]{Lemma}

\newtheorem{exam}[thm]{Example}
\newtheorem{rem}[thm]{Remark}

\def\proof{{\parindent0pt {\bf Proof.\ }}}

\theoremstyle{definition}

\theoremstyle{remark}

\theoremstyle{Definition and Notation}

\begin{document}

\bibliographystyle{amsplain}


\title[n-coherence and (n, d)- properties in amalgamated algebra]{$\textbf{n}$-coherence and $\textbf{(n, d)}$-properties in amalgamated algebra along an ideal}

\author{Karima Alaoui Ismaili}
\address{Karima Alaoui Ismaili\\Department of Mathematics, Faculty of Science and Technology of Fez, Box 2202,
University S.M. Ben Abdellah Fez, Morocco.
$$ E-mail\ address:\ alaouikarima2012@hotmail.fr$$}

\author{Najib Mahdou}
\address{Najib Mahdou\\Department of Mathematics, Faculty of Science and Technology of Fez, Box 2202,
University S.M. Ben Abdellah Fez, Morocco.
$$E-mail\ address:\ mahdou@hotmail.com$$}


\begin{abstract}
Let  $f: A\rightarrow B$ be a ring homomorphism and let $J$ be an ideal of $B$. The purpose of this article is to examine the transfer of the properties of $n$-coherence and strong $n$-coherence from a ring $A$ to his amalgamated algebra $A\bowtie^{f} J$. Also, we investigate the $(n,d)$-property of the amalgamated algebra $A \bowtie^{f}J$, to resolve Costa's first conjecture.
\end{abstract}

\maketitle

 \begin{section} {Introduction}
 Throughout this paper, all rings are commutative with identity element, and all modules are unitary.\par
Let $A$ and $B$ be two rings, let $J$ be an ideal of $B$ and let $f: A\rightarrow B$ be a ring homomorphism. In this setting, we can consider the following subring of $A \times B$:\\
\begin{center}
$A\bowtie^{f}J = \{(a,f(a)+j) \diagup a \in A, j \in J\}$
\end{center}
called the amalgamation of $A$ with $B$ along $J$ with respect to $f$ (introduced and studied by D'Anna, Finocchiaro, and Fontana in \cite{AFF1,AFF2}). This construction is a generalization of the amalgamated duplication of a ring along an ideal (introduced and studied by D'Anna and Fontana in \cite{A, AF, AF2} and denoted by $A \bowtie I$). Moreover, other classical constructions (such as the $A+XB[X]$, $A+XB[[X]]$, and the $D+M$ constructions) can be
studied as particular cases of the amalgamation \cite[Examples 2.5 \& 2.6]{AFF1} and other classical constructions, such as the Nagata's idealization and the CPI extensions (in the sense of
Boisen and Sheldon \cite{Boi}) are strictly related to it (see \cite[Example 2.7 \& Remark 2.8]{AFF1}).\par
Let $R$ be a commutative ring. For a nonnegative integer $n$, an $R$-module $E$ is called $n$-presented
if there is an exact sequence of $R$-modules: \\
$$F_n \rightarrow  F_{n-1} \rightarrow \ldots F_1  \rightarrow F_0  \rightarrow E \rightarrow 0$$ \\
where each $F_i$ is a finitely generated free $R$-module. In particular, $0$-presented and $1$-presented $R$-modules are,
respectively, finitely generated and finitely presented $R$-modules.\par
The ring $R$ is $n$-coherent if each
$(n-1)$-presented ideal of $R$ is $n$-presented, and $R$ is strong $n$-coherent ring if each $n$-presented
$R$-module is $(n+1)$-presented \cite{DKM, DKMS} (This terminology is not the same as that of
Costa $(1994)$ \cite{Co}, where Costa's $n$-coherence is our strong $n$-coherence). In particular,
$1$-coherence'' coincides with coherence, and one may view $0$-coherence as Noetherianity.
Any strong $n$-coherent ring is $n$-coherent, and the converse holds for $n =1$ or for
coherent rings \cite[Proposition 3.3]{DKMS}. \par
In 1994, Costa \cite{Co} introduced a doubly filtered set of classes of rings in order to categorize the structure of non-Noetherian rings: for non-negative integers $n$ and $d$, we say that a ring $R$ is an
$(n, d)$-ring if $pd_{R}(E) \leq d $ for each $n$-presented $R$-module $E$ (as usual, $pd_{R}(E)$ denotes the projective dimension of $E$ as an $R$-module). An integral domain with this property will be called
an $(n, d)$-domain. For example, the $(n, 0)$-domains are the fields, the
$(0, 1)$-domains are the Dedekind domains, and the $(1, 1)$)-domains are the
Pr\"ufer domains \cite{Co}. Every $(n,d)$-ring is strong $(sup\{n,d\})$-coherent and every $(n,d)$-domain is strong $(sup\{n,d-1\})$-coherent \cite[Theorem 2.2]{Co}.\par
We call a commutative ring an $n$-Von Neumann regular ring if it is an $(n, 0)$-ring. Thus, the $1$-Von Neumann regular rings are the Von Neumann regular rings \cite[Theorem 1.3]{Co}.\par
In \cite{Co}, Costa asks whether there is an $(n, d)$-ring which is neither an $(n, d-1)$-ring nor an $(n-1, d)$-ring for each integers $n, d \geq 0$. The answer is affirmative for $(0, d)$-ring and $(1, d)$-ring for each integer $d \geq 0$ (\cite{Co}). \par
Again in \cite{Co}, Costa gives examples of $(2,1)$-domains which are neither $(2, 0)$-domains (fields) nor $(1, 1)$-domains (Pr\"ufer), and in \cite{CK}, Costa and Kabbaj give examples of $(2, 2)$-domains which are neither $(2, 1)$-domains nor $(1, 2)$-domains. Later, in \cite{M1, M4}, the author gives a class of $(2, d)$-domains which are neither $(2, d-1)$-domains nor $(1, d)$-domains for each integer $d \geq 1$, and a class of $(2, d)$-rings (not domains) which are neither $(2, d-1)$-rings (for $d \geq 1$) nor $(1, d)$-rings for each integer $d \geq 0$. Next, in \cite{KM1}, the authors construct a class of $(3, d)$-rings which are
neither $(3, d-1)$-rings (for $d \geq 1$) nor $(2, d)$-rings for each integer $d \geq 0$. Finally, in \cite{M3}, the author gives a sufficient condition to resolve Costa's first conjecture for each positive integer $n$ and $d$ with $n\geq 4$. The second main goal of this paper is the constructions of the second class of $(3, d)$-rings for each integer $d \geq 0$, and $(2, d)$-rings for $d \leq 2$, after the first class of cost's conjecture given by the authors in \cite{KM1, M1, M4}. \par
Let $A$ be a ring, $E$ be an $A$-module, and $R:=A\propto E$ be the set of pairs $(a,e)$ with pairwise addition and multiplication given by $(a,e) (b,f)= (ab,af + be)$. $R$ is called the trivial ring extension of $A$ by $E$ (also called the idealization of $E$ over $A$). Considerable work, part of it summarized in Glaz \cite{Gz2} and Huckaba \cite{H}, has been concerned with trivial ring extensions. These have proven to be useful in solving many open problems and conjectures for various contexts in (commutative and noncommutative) ring theory. See for instance \cite{BKM, BG, Gz2, H, KM1, KM2, PR}. \par
The first section of this work examines the transfer of the properties of $n$-coherence and strong $n$-coherence to the amalgamated algebra $A\bowtie^{f}J$. Thereby, new examples are provided which, particularly, enriches the current literature with new classes of $n$-coherent rings ($n \geq 2$) that are non-coherent rings.\par

\bigskip
\bigskip
 \begin{section}{Transfer of the properties of (strong) $n$-coherence ($n \geq 1$)}
\bigskip


\bigskip
The first main result of this section (Theorem \ref{thm1}) examines the transfer of the properties of strong $n$-coherence and $n$-coherence ($n \geq 1$) to the amalgamation algebra along an ideal $A\bowtie^{f}J$ issued from local rings.\par
First, it is worthwhile recalling that the function $f^{n}: A^{n}\rightarrow B^{n}$ defined by $f^{n}((\alpha_{i})_{i=1}^{i=n}) = (f(\alpha_{i}))_{i=1}^{i=n}$
is a ring homomorphism, $(A\bowtie^{f}J)^{n} \cong A^{n}\bowtie^{f^{n}}J^{n}$ and $f^{n}(\alpha a)= f(\alpha)f^{n}(a)$ for all $\alpha \in A$ and $a \in A^{n}$
(see \cite{AM1}).\\

\bigskip
Next, before we announce the main result of this section (Theorem
\ref{thm1}), we make the following useful remark. \\

\begin{rem} \label{rem1}
Let $(A, M)$ be a local ring, $f: A\rightarrow B$ be a ring homomorphism, and let $J$ be a proper ideal of $B$
 such that $J^{2}=0$. Then $A\bowtie^{f}J$ is a local ring and $M\bowtie^{f}J$ is its maximal ideal.
 \end{rem}

Indeed, by \cite[Proposition 2.6 (5)]{AFF2}, $Max(A\bowtie^{f}J)=\{m\bowtie^{f}J \diagup m \in Max(A)\} \cup \{\overline{Q}\}$ with $Q \in Max(B)$ not containing $V(J)$ and $\overline{Q}:= \{(a,f(a)+j) \diagup a\in A, j\in J, f(a)+j\in Q \}$.
 Since $J^{2}=0$, then $J \subseteq Rad(B)$, and then $J \subseteq Q$ for all $Q \in Max(B)$.
  So, $Max(A\bowtie^{f}J)=\{m\bowtie^{f}J \diagup m \in Max(A)\}=M\bowtie^{f}J$
  since $(A, M)$ is a local ring. Therefore $(A\bowtie^{f}J, M\bowtie^{f}J)$ is a local ring, as desired. \\
  \qed
\bigskip
\begin{thm} \label{thm1}
Let $(A, M)$ be a local ring, $f: A\rightarrow B$ be a ring
homomorphism, and let $J$ be a proper ideal of $B$ such that $J$
is a finitely generated ideal of $(f(A)+J)$, $J^{2}=0$ and
$f(M)J=0$. Then:
\begin{enumerate}
\item $A\bowtie^{f}J $ is strong $n$-coherent ring if and only if so is $A$.
\item $A\bowtie^{f}J $ is $n$-coherent ring if and only if so is $A$.
\end{enumerate}
\end{thm}
\bigskip
The proof of Theorem \ref{thm1} draws on the following results. \\

\begin{lem} \label{lem1}
Let $(A, M)$ be a local ring, $f: A\rightarrow B$ be a ring homomorphism, and let $J$ be a proper ideal of $B$ such that $J^{2}=0$ and $f(M)J=0$. Let $p, n \in \N$, and let $U$ be a sub-module of $M^{p}$. Then $U \bowtie^{f^{p}}J^{p}$ is a $n$-finitely presented $(A\bowtie^{f}J )$-module if and only if $U$ is a $n$-finitely presented $A$-module, and $J$ is a finitely generated ideal of $f(A)+ J$.
\end{lem}

\proof Proceed by induction on $n$. The property is true for $n = 0$. Indeed, by \cite[Lemma 2.4]{AM1}, it remains to show that
if $U \bowtie^{f^{p}}J^{p} := \sum_{i=1}^{i=r} (A\bowtie^{f}J)(u_{i}, f^{p}(u_{i})+ k_{i})$ is a finitely generated ideal of
$A\bowtie^{f}J$, where, $u_{i} \in U$ and $k_{i} \in J^{p}$ for all $ 1 \leq i \leq  r$, then $J$ is a finitely generated ideal
of $(f(A)+J)$. Clearly, $\sum_{i=1}^{i=r}(f(A)+J)k_{i} \subseteq J^{p}$. Let $k \in J^{p}$. Then,
$(0,k) = \sum_{i=1}^{i=r} (\alpha_{i},f(\alpha_{i}) + j_{i})(u_{i}, f^{p}(u_{i})+ k_{i})$ for some $\alpha_{i} \in A$ and $j_{i}\in J$.
So $\sum_{i=1}^{i=r}\alpha_{i}u_{i}=0$, and $k = \sum_{i=1}^{i=r}(f(\alpha_{i}) + j_{i})(f^{p}(u_{i})+ k_{i})= \sum_{i=1}^{i=r}(f(\alpha_{i}) + j_{i})k_{i} \in \sum_{i=1}^{i=r}(f(A)+J)k_{i}$ (since $f^{p}(U)J \subseteq f^{p}(M^{p})J=0$). Thus  $J^{p} \subseteq \sum_{i=1}^{i=r}(f(A)+J)k_{i}$. Therefore, $J^{p}=\sum_{i=1}^{i=r}(f(A)+J)k_{i} $, and then $J$ is a finitely generated ideal of $f(A)+J$. Assume that the property is true for $n$, and assume that $U \bowtie^{f^{p}}J^{p}= \sum_{i=1}^{i=r}A\bowtie^{f}J(u_{i},f^{p}(u_{i})+ k_{i})$ is a $(n+1)$-finitely presented $(A\bowtie^{f}J )$-module, where $u_{i} \in U$ and $k_{i} \in J^{p}$ for all $ i\in\{1,.....r\}$. Clearly, $U=\sum_{i=1}^{i=r}Au_{i}$. We may assume that $\{(u_{i},f^{p}(u_{i})+ k_{i})_{i=1}^{i=r} \}$ is a minimal generating set of $U \bowtie^{f^{p}}J^{p}$. Consider the exact sequence of $A$-modules:
\begin{center}
$ 0\rightarrow Kerv \rightarrow A^{r} \rightarrow U \rightarrow 0 $ \hspace{2 cm} (1)
\end{center}
where $v((\alpha_{i})_{i=1}^{i=r}) = \sum_{i=1}^{i=r}\alpha_{i}u_{i}$. On the other hand consider the exact sequence of ($A\bowtie^{f}J$)-modules:
\begin{center}
$ 0\rightarrow Keru \rightarrow (A\bowtie^{f}J)^{r} \rightarrow U \bowtie^{f^{p}}J^{p} \rightarrow 0 $ \hspace{2 cm} (2)
\end{center}
where
\begin{eqnarray}
\nonumber  u((\alpha_{i},f(\alpha_{i}) + j_{i})_{i=1}^{i=r}) &=& \sum_{i=1}^{i=r}(\alpha_{i},f(\alpha_{i}) + j_{i})(u_{i},f^{p}(u_{i})+ k_{i})\\
\nonumber    &=& (\sum_{i=1}^{i=r}\alpha_{i}u_{i},\sum_{i=1}^{i=r}f(\alpha_{i})(f^{p}(u_{i})+ k_{i})).
\end{eqnarray}
 Then $Keru = \{(\alpha_{i},f(\alpha_{i}) + j_{i})_{i=1}^{i=r} \in (A\bowtie^{f}J)^{r}\diagup \sum_{i=1}^{i=r}\alpha_{i}u_{i}=0, \sum_{i=1}^{i=r}f(\alpha_{i})(f^{p}(u_{i})+ k_{i}) = 0\}$.\\
 So, $ Keru= \{((\alpha_{i})_{i=1}^{i=r},f^{r}((\alpha_{i})_{i=1}^{i=r}) + (j_{i})_{i=1}^{i=r}) \in A^{r}\bowtie^{f^{r}}J^{r}\diagup \sum_{i=1}^{i=r}\alpha_{i}u_{i}=0, \sum_{i=1}^{i=r}f(\alpha_{i})k_{i}=0 \}$. Since $A\bowtie^{f}J$ is a local ring by Remark \ref{rem1} and $\{(u_{i},f^{p}(u_{i})+ k_{i})_{i=1}^{i=r} \}$ is a minimal generating set of $U \bowtie^{f^{p}}J^{p}$, then $Keru \subseteq M^{r}\bowtie^{f^{r}}J^{r}$. So,
  $ Keru= \{((\alpha_{i})_{i=1}^{i=r},f^{r}((\alpha_{i})_{i=1}^{i=r}) + (j_{i})_{i=1}^{i=r}) \in A^{r}\bowtie^{f^{r}}J^{r}\diagup (\alpha_{i})_{i=1}^{i=r} \in Kerv\}$.
Therefore $Keru = Kerv \bowtie^{f^{r}}J^{r}$.\\
Since $U \bowtie^{f^{p}}J^{p}$ is a $(n+1)$-finitely presented $(A\bowtie^{f}J)$-module, then $Keru$ is a $n$-finitely presented $(A\bowtie^{f}J )$-module (by a sequence (2)). So, $Kerv$ is a $n$-finitely presented $A$-module and $J$ is a finitely generated ideal of $(f(A)+J)$ by induction (since $Kerv \subseteq M^{r}$). Thus, $U$ is a $(n+1)$-finitely presented $A$-module (by a sequence (1)). Conversely, assume that $U$ is a $(n+1)$-finitely presented $A$-module and $J$ is a finitely generated ideal of $(f(A)+J)$, then $U \bowtie^{f^{p}}J^{p}$ is a finitely generated $(A\bowtie^{f}J)$-module by induction, and then $U \bowtie^{f^{p}}J^{p}= \sum_{i=1}^{i=r}A\bowtie^{f}J(u_{i},f^{p}(u_{i})+ k_{i})$, where $u_{i} \in U$ and $k_{i} \in J^{p}$ for all $ i\in\{1,.....r\}$. It is obvious that $U=\sum_{i=1}^{i=r}Au_{i}$.
Since $U$ is a $(n+1)$-finitely presented $A$-module, then $Kerv$ is a $n$-finitely presented $A$-module (by a sequence (1)). So, $Keru (=Kerv \bowtie^{f^{r}}J^{r})$ is a $n$-finitely generated $(A\bowtie^{f}J)$-module by induction, and then $U \bowtie^{f^{p}}J^{p}$ is a $(n+1)$-finitely presented ($A\bowtie^{f}J $)-module (by a sequence (2)), as desired. \\
 \qed
\bigskip
\begin{lem} \label{lem2}
Let $(A, M)$ be a local ring, $f: A\rightarrow B$ be a ring homomorphism, and let $J$ be a proper ideal of $B$ such that $J^{2}=0$ and $f(M)J=0$. Let $p, n \in \N^{\ast}$, $W := \sum_{i=1}^{i=r} (A\bowtie^{f}J )(u_{i},f^{p}(u_{i})+ k_{i})$ where, $u_{i} \in M^{p}$, $k_{i} \in J^{p}$, and let $U=\sum_{i=1}^{i=r}Au_{i}$. Then $W$ is a $n$-finitely presented $(A\bowtie^{f}J )$-module if and only if $U$ is a $n$-finitely presented $A$-module and $J$ is a finitely generated ideal of $f(A)+ J$.
\end{lem}

\proof   Let $W := \sum_{i=1}^{i=r} (A\bowtie^{f}J )(u_{i},f^{p}(u_{i})+ k_{i})$ and $U=\sum_{i=1}^{i=r}Au_{i}$. We may assume that  $\{(u_{i},f^{p}(u_{i})+ k_{i})_{i=1}^{i=r} \}$ is a minimal generating set of $W$. Consider the exact sequence of $A$-modules:
\begin{center}
$ 0\rightarrow Kerv \rightarrow A^{r} \rightarrow U \rightarrow 0 $ \hspace{2 cm} (1)
\end{center}
where $v((\alpha_{i})_{i=1}^{i=r}) = \sum_{i=1}^{i=r}\alpha_{i}u_{i}$. On the other hand consider the exact sequence of ($A\bowtie^{f}J$)-modules:
\begin{center}
$ 0\rightarrow Keru \rightarrow (A\bowtie^{f}J)^{r} \rightarrow W \rightarrow 0 $ \hspace{2 cm} (2)
\end{center}
where
\begin{eqnarray}
\nonumber u((\alpha_{i},f(\alpha_{i}) + j_{i})_{i=1}^{i=r}) &=&  \sum_{i=1}^{i=r}(\alpha_{i},f(\alpha_{i}) + j_{i})(u_{i},f^{p}(u_{i})+ k_{i})\\
\nonumber   &=& (\sum_{i=1}^{i=r}\alpha_{i}u_{i}, \sum_{i=1}^{i=r}f(\alpha_{i})(f^{p}(u_{i})+ k_{i}))
\end{eqnarray}
Then $Keru = \{(\alpha_{i},f(\alpha_{i}) + j_{i})_{i=1}^{i=r} \in (A\bowtie^{f}J)^{p}\diagup \sum_{i=1}^{i=r}\alpha_{i}u_{i} =0,\sum_{i=1}^{i=r}f(\alpha_{i})k_{i}= 0\}$. Since $Keru \subseteq M^{r}\bowtie^{f^{r}}J^{r}$, then $Keru = Kerv \bowtie^{f^{r}}J^{r}$.\\
By Lemma \ref{lem1}, $Keru$ is a $n$-finitely presented $(A\bowtie^{f}J)$-module if and only if $Kerv$ is a $n$-finitely presented $A$-module and $J$ is a finitely generated ideal of $f(A) + J$ (since $Kerv \subseteq M^{r}$). So, $W$ is a $(n+1)$-finitely presented $(A\bowtie^{f}J )$-module if and only if $U$ is a $(n+1)$-finitely presented $A$-module and $J$ is a finitely generated ideal of $f(A)+ J$, as desired. \\
\qed

\bigskip
\noindent {\bf Proof of Theorem \ref{thm1}.} \par {\bf (1)} Recall
that $R$ is a strong $n$-coherent ring if and only if every
$(n-1)$-presented submodule of a finitely generated free
$R$-module is $n$-presented. Assume that $A\bowtie^{f}J$ is strong
$n$-coherent ring and let $U:= \sum_{i=1}^{i=r} Au_{i}$ be a
$(n-1)$-finitely presented $A$-module, where $u_{i} \in M^{p}$,
then $W:\sum_{i=1}^{i=r} (A\bowtie^{f}J )(u_{i},f^{p}(u_{i}))$ is
a $(n-1)$-finitely
  presented $(A\bowtie^{f}J)$-module by Lemma \ref{lem2}. So, $W$ is $n$-finitely presented $(A\bowtie^{f}J)$-module since $A\bowtie^{f}J$ is strong $n$-coherent ring. Therefore $U$ is $n$-finitely presented $A$-module by Lemma \ref{lem2}.
  Thus, $A$ is strong $n$-coherent ring. Conversely, assume that $A$ is strong $n$-coherent ring and
  let $W:= \sum_{i=1}^{i=r} A\bowtie^{f}J(u_{i}, f^{p}(u_{i})+k_{i})$ be a $(n-1)$-finitely presented ($A\bowtie^{f}J$)-module, where $u_{i} \in M^{p}$, and $k_{i} \in J^{p}$.
   Then $U:= \sum_{i=1}^{i=r} Au_{i}$ is a $(n-1)$-finitely presented $A$-module by Lemma \ref{lem2}. So, $U$ is $n$-finitely presented $A$-module since $A$ is strong $n$-coherent ring,
   and then $W$ is a $n$-finitely presented $(A\bowtie^{f}J)$-module by Lemma \ref{lem2}. Hence, $A\bowtie^{f}J$ is strong $n$-coherent ring. \\
  {\bf (2)} The same reasoning as in the proof of (1) shows that $A\bowtie^{f}J $ is $n$-coherent ring if and only if so is $A$.\\
\qed

 \bigskip
 The following corollaries are an immediate consequence of Theorem
 \ref{thm1}. \\

 \begin{cor} \label{cor1}
 Let $(A, M)$ be a local ring, $f: A\rightarrow B$ be a ring homomorphism, and let $J$ be a proper ideal of $B$ such that $J$ is a finitely
 generated ideal of $(f(A)+J)$, $J^{2}=0$ and $f(M)J=0$. Then, $A\bowtie^{f}J$ is a (strong) $n$-coherent ring which is
  non-(strong) $(n-1)$-coherent if and only if so is $A$.
 \end{cor}
 \bigskip
 \begin{cor} \label{cor2}
 Let $(A,M)$ be a local ring and $I$ be a finitely generated ideal of $A$ such that $MI=0$. Then:
\begin{enumerate}
\item $A\bowtie I $ is a strong $n$-coherent ring if and only if so is $A$.
\item $A\bowtie I $ is a $n$-coherent ring if and only if so is $A$.
\end{enumerate}
 \end{cor}
\bigskip
In particular, Theorem \ref{thm1} enriches the literature with new
examples of $2$-coherent rings which are non-coherent rings. \\

\begin{exam}
Let $k$ be a field, $K$ be a field containing $k$ as a subfield such that $[K:k]= \infty$, $T=K[[X]]=K+M$, where $X$ is an indeterminate over $K$ and $M=XT$ be the maximal ideal of $T$, $A=k+M$. \\
Set $B:= A/M^{2}$, $J:=Bf(m)$ be an ideal of $B$, where $m \in M$ such that $f(m)\neq 0$, and consider the canonical ring homomorphism $f: A\rightarrow B$ ($f(x) = \overline{x}$). Then:
\begin{enumerate}
\item By Theorem \ref{thm1}, $A\bowtie^{f}J$ is $2$-coherent ring since $A$ is by \cite[Corollary 5.2]{Co}.
\item By \cite[Theorem 4.1.5]{Gz2}, $A\bowtie^{f}J$ is non-coherent ring since $A$ is non-coherent domain by \cite[Corollary 5.2]{Co}.
\end{enumerate}
\end{exam}
\bigskip
\begin{exam}
Let $K$ be a field and $E$ be a $K$-vector space of infinite dimension. Let $A:= K\propto E$ be trivial extension ring of $K$ by $E$. Set $I:= 0 \propto E^{'}$, where $E^{'}$ is a finitely generated $K$-subspace of $E$. Then:
\begin{enumerate}
\item By Corollary \ref{cor2}, $A\bowtie I $ is a $2$-coherent ring since $A$ is by \cite[Theorem 3.4]{M1}.
\item By \cite[Theorem 4.1.5]{Gz2} $A\bowtie I $ is non-coherent ring since $A$ is non-coherent ring by \cite[Theorem 2.6]{KM2}.
\end{enumerate}
\end{exam}

\end{section}
\bigskip
 \begin{section}{Transfer of the $(n,d)$-property }
\bigskip
For integers $n, d \geq 0$, Costa asks in \cite{Co} whether there
is an $(n, d)$-ring which is neither an $(n,d-1)$-ring nor an
$(n-1,d)$-ring? The answer is affirmative for $(0, d)$-rings, $(1,
d)$-rings, $(2, d)$-rings, $(3, d)$-rings, and $(n+4,d)$-rings for
all integers $n, d$ (See for instance \cite{Co, CK, KM1, KM2, M1,
M4, Z, M3}). The goal of This section is to give a second class of
cost's conjecture after the first given by the authors in
\cite{KM1, M1, M4}. At the end of this work, we will be able to
give examples of $(2, d)$-rings which is neither a $(1,d)$-ring
($d=0, 1, 2$) nor a $(2,d-1)$-ring ($d=1, 2$), and examples of
$(3, d)$-rings which is neither a $(2,d)$-ring ($d\geq 0$) nor a
$(3,d-1)$-ring ( for each integer $d$). \\
The following Theorem \ref{thm2} allows us to provide new examples of $(2,0)$, $(2,1)$, and $(2,2)$-ring, to resolve cost's conjecture.\\
\bigskip
\begin{thm} \label{thm2}
Let $(A, M)$ be a local ring, $f: A\rightarrow B$ be a ring homomorphism, and let $J$ be a proper ideal of $B$ such that $f(M)J=0$ and $J^{2}=0$. Then:
\begin{enumerate}
\item
\begin{enumerate}
\item $A\bowtie^{f}J $ is a $(2,0)$-ring provided $J$ is a not finitely generated ideal of $f(A)+ J$.
\item If $A\bowtie^{f}J $ is a $(2,0)$-ring, then $M$ is not a finitely generated ideal of $A$ or $J$ is not a finitely generated ideal of $f(A)+ J$.
\item Assume that $M$ is a finitely generated ideal of $A$. Then $A\bowtie^{f}J $ is a $(2,0)$-ring if and only if $J$ is a not finitely generated ideal of $f(A)+ J$.
\end{enumerate}
\item $A\bowtie^{f}J $ is a non-$(1,2)$-ring. In particular, $A\bowtie^{f}J $ is a non-von
Neumann regular ring.
\end{enumerate}
\end{thm}

\bigskip
\noindent {\bf Proof of Theorem \ref{thm2}} \par
{\bf (1) (a)} Let $K$ be a $2$-finitely presented $(A\bowtie^{f}J)$-module, and let $\{(k_{i})_{i=1}^{i=p} \}$ be a minimal generating set of $K$. We want to show that $K$ is a projective ($A\bowtie^{f}J$)-module. For this, consider the exact sequence of ($A\bowtie^{f}J$)-modules:
\begin{center}
$ 0\rightarrow Keru (=H) \rightarrow (A\bowtie^{f}J)^{p} \rightarrow  K \rightarrow 0 $
\end{center}
where $u((\alpha_{i})_{i=1,.....,p}) = \sum_{i=1}^{i=p}\alpha_{i}k_{i}$.
Then $H \subset (M\bowtie^{f}J)^{p}$ since $\{(k_{i})_{i=1}^{i=p} \}$ is a minimal generating set of $K$, $A\bowtie^{f}J$ is a local ring with maximal ideal $M\bowtie^{f}J$ by Remark \ref{rem1}. We prove that $H=0$. Otherwise, $H \neq 0$. Let $(m_{i},f^{p}(m_{i})+k_{i})_{\{i=1,...,r\}}$ be a minimal generating set of H, where $m_{i} \in M^{p}, k_{i} \in J^{p}$ for each $i=1,.......,r$. Consider the exact sequence of $(A\bowtie^{f}J)$-modules
\begin{center}
$ 0\rightarrow Kerv \rightarrow (A\bowtie^{f}J)^{r} \rightarrow H  \rightarrow 0 $
\end{center}
where $v((\alpha_{i},f(\alpha_{i})+j_{i})_{i=1,.....,r}) = \sum_{i=1}^{i=r}(\alpha_{i}m_{i},(f(\alpha_{i})+j_{i})(f^{p}(m_{i})+k_{i}))=
(\sum_{i=1}^{i=r}\alpha_{i}m_{i},\sum_{i=1}^{i=r}f(\alpha_{i})(f^{p}(m_{i})+k_{i}))$. But $Kerv \subset (M\bowtie^{f}J)^{r}$. Hence
$Kerv = U\bowtie^{f^{r}}J^{r}$. Where $U=\{(\alpha_{i})_{i=1,...,r} \in A^{r} \diagup \sum_{i=1}^{i=r}\alpha_{i}m_{i}=0 \}$. Since
$K$ is a $2$-presented $(A\bowtie^{f}J)$-module, then $Kerv$ is a finitely generated $(A\bowtie^{f}J)$-module. So, $J$ is a finitely
generated ideal of $(f(A)+ J)$ by Lemma \ref{lem1} (since $U \subseteq M^{r}$). A contradiction since $J$ is a not finitely generated
ideal of $f(A)+ J$ by hypothesis. So $H = 0$. Hence, $K \cong (A\bowtie^{f}J)^{p}$ is a projective $(A\bowtie^{f}J)$-module.\\

{\bf (b)} Let $k (\neq 0) \in J$, $I=(A\bowtie^{f}J)(0,k)$, and consider the exact sequence of ($A\bowtie^{f}J$)-modules:
\begin{center}
$ 0\rightarrow Keru  \rightarrow A\bowtie^{f}J \rightarrow  I \rightarrow 0 $
\end{center}
where $u(\alpha,f(\alpha)+j) = (0, f(\alpha)k)$. So, $Keru=M\bowtie^{f}J$ is a not finitely generated ideal of $A\bowtie^{f}J$ since $A\bowtie^{f}J$ is a
$(2,0)$-ring, and $I$ is not projective ideal of $A\bowtie^{f}J$ (otherwise, $I$ is a free ideal of $A\bowtie^{f}J$ since $A\bowtie^{f}J$ is a local
ring, absurd since $(0,k)I=0$). Therefore, $M$ is a not finitely generated ideal of $A$ or $J$ is a not finitely generated ideal of $f(A)+J$ by
Lemma \ref{lem1}.\\

{\bf (c)} Follows immediately from (a) an (b).\\

{\bf (2)} Let $I=(A\bowtie^{f}J)(0,k)$, where $k (\neq 0) \in J$, and consider the exact sequence of ($A\bowtie^{f}J$)-modules:
\begin{center}
$ 0\rightarrow M\bowtie^{f}J  \rightarrow A\bowtie^{f}J \rightarrow  I \rightarrow 0 $
\end{center}
Since $A\bowtie^{f}J$ is a local ring and $M\bowtie^{f}J$ is not a free ideal of $A\bowtie^{f}J$ (since $(0,k) (M\bowtie^{f}J)=0$), then $M\bowtie^{f}J$ is
not a projective ideal of $A\bowtie^{f}J$. So, $pd_{A\bowtie^{f}J}(\frac{A\bowtie^{f}J}{I}) > 2$ (i.e. $A\bowtie^{f}J$ is a not a $(1,2)$-ring), and
this completes the proof of Theorem \ref{thm2}. \\
\qed

\bigskip
The next corollary is an immediate consequence of Theorem
\ref{thm2}. \\

\begin{cor} \label{cor3}
Let $(A,M)$ be a local ring, $f: A\rightarrow B$ be a ring homomorphism, and let $J$ be a proper ideal of $B$ such that $f(M)J=0$ and $J^{2}=0$.
Let $A_{1}=A\bowtie^{f}J$, $d\leq 2$ be an integer, $A_{2}$ be a Noetherian ring of global dimension $d$, and let $C=A_{1} \times A_{2}$ the direct
product of $A_{1}$ and $A_{2}$.\\ Assume that $J$ is a not finitely generated ideal of $f(A)+J$. Then $C$ is a $(2,d)$-ring which is neither a
$(1,d)$-ring $(d=0, 1, 2)$ nor a $(2,d-1)$-ring ($d=1, 2$).
\end{cor}
\proof  Let $d\leq 2$ be an integer. By \cite[Theorem 2.4]{M1},  $C=A_{1} \times A_{2}$ is a $(2,d)$-ring since $A_{i}$ is a $(2,d)$-ring for each
$i=1, 2$, and $C$ is not $(1,d)$-ring since $A_{1}$ is not $(1,d)$-ring by Theorem \ref{thm2} (2). It remains to show that $C$ is not a $(2, d-1)$-ring
for $1\leq d\leq 2$. Assume that $C$ is a $(2,d-1)$-ring. By \cite[theorem 2.4]{M1}, $A_{2}$ is also a $(2,d-1)$-ring and then is a $(0,d-1)$-ring
by \cite[Theorem 2.4]{Co} since $A_{2}$ is a Noetherian ring. Thus, $gldim(A_{2}) \leq d-1$, but this is a contradiction since $gldim(A_{2})=d$, as desired.\\
\qed

\bigskip

Now, we are able to give a new examples of $(2,0)$, $(2,1)$, and
$(2,2)$-ring, to resolve cost's conjecture. \\

\begin{exam}
Let $(R,m)$ be a local ring such that $m$ is a finitely generated
ideal of $R$ and $m^{2}=0$ (for instance
($(R,m)=(\frac{\mathbb{Z}}{4\mathbb{Z}},
\frac{2\mathbb{Z}}{4\mathbb{Z}}$)), $E$ be an $\frac{R}{m}$-vector
space with finite rank, and let $A:=R\propto E$ and $M:=m \propto
E$. Let $E^{'}$ be a $\frac{A}{M}$-vector space of infinite
dimension, $B:=A\propto E^{'}$, $J:=0\propto E^{'}$, and consider
the ring homomorphism $f: A \rightarrow B$ $(f(a)=(a,0))$. Let $K$
be a field, $A_{1}=A\bowtie^{f}J$, and let  $A_{2}:=K[X_{1}]$ and
$A_{3}:=K[X_{1}, X_{2}]$ where $X_{1}, X_{2}$ are indeterminate
over $K$. Then, by Corollary \ref{cor3}:
\begin{enumerate}
\item $A_{1}$ is a $(2,0)$ ring that is not a $(1,0)$-ring.
\item $A_{1} \times A_{2}$ is a $(2,1)$-ring which is neither a $(1,1)$-ring nor a $(2,0)$-ring.
\item $A_{1} \times A_{3}$ is a $(2,2)$-ring which is neither a $(1,2)$-ring nor a $(2,1)$-ring.
\end{enumerate}
\end{exam}

\bigskip

 The aim of Theorem \ref{thm3} is to construct a class of $(3, d)$-rings which are neither $(3, d-1)$-rings (for each positive integer $d$) nor
 $(2, d)$-rings for each integer $d \geq 0$. \\

\begin{thm} \label{thm3}
Let $(A, M)$ be a local ring, $f: A\rightarrow B$ be a ring homomorphism, and let $J$ be a proper ideal of $B$ such that $f(M)J=0$ and $J^{2}=0$.
\begin{enumerate}
\item $A\bowtie^{f}J $ is a $(3,0)$-ring provided $M$ is not a finitely generated ideal of $A$.
\item Assume that $J$ is a finitely generated ideal of $f(A)+ J$ and $M$ contains a regular element. Then:
\begin{enumerate}
\item $A\bowtie^{f}J $ is not a $(2,d)$-ring for each positive integer $d$.
\item Let $A_{1}=A\bowtie^{f}J$, $d$ be an integer, $A_{2}$ be a Noetherian ring of global dimension $d$, and let $C=A_{1} \times A_{2}$ the direct product of $A_{1}$ and $A_{2}$. Assume that $M$ is not finitely generated ideal of $A$. Then $C$ is a $(3,d)$-ring which is neither a $(2,d)$-ring $(d\geq 0)$ nor a $(3,d-1)$-ring $(d\geq 1)$.
\end{enumerate}
\end{enumerate}
\end{thm}
\bigskip
The proof of this theorem requires the next result. \\

\begin{lem} \label{lem3}
Let $(A, M)$ be a local ring, $f: A\rightarrow B$ be a ring homomorphism, and let $J$ be a proper ideal of $B$ such that $J$ is a finitely generated
ideal of $f(A)+ J$, $f(M)J=0$ and $J^{2}=0$. Then $pd_{A\bowtie^{f}J}(M\bowtie^{f}J)$ and $pd_{A\bowtie^{f}J}(\{0\} \times J)$ are infinite.
\end{lem}
\proof  Consider the exact sequence of $(A\bowtie^{f}J)$-modules
\begin{center}
$ 0\rightarrow M\bowtie^{f}J \rightarrow A\bowtie^{f}J \rightarrow \frac{A\bowtie^{f}J}{M\bowtie^{f}J} \rightarrow 0 $
\end{center}
We claim that $\frac{A\bowtie^{f}J}{M\bowtie^{f}J}$ is not projective. Otherwise, the sequence splits.
Hence, $M\bowtie^{f}J$ is generated by an idempotent element $(m, f(m)+k) = (m, f(m)+k)(m, f(m)+k) = (m^{2}, f(m^{2}))$. So $M\bowtie^{f}J = (A\bowtie^{f}J)(m,f(m)) = Am \bowtie^{f}0$, the desired contradiction (since $J \neq 0$). It follows from the above sequence that
\begin{center}
$pd_{A\bowtie^{f}J}(M\bowtie^{f}J)+1=pd_{A\bowtie^{f}J}(\frac{A\bowtie^{f}J}{M\bowtie^{f}J})$ \hspace{2 cm} (1)
\end{center}
Let $\{(m_{i})_{i \in I}\}$ be a set of generators of $M$, where $m_{i} \in M$ for all $i \in I$, and let $\{(g_{i})_{i=1}^{i=p}\}$ be a minimal generating set of $J$. Consider the exact sequence of $(A\bowtie^{f}J)$-modules
\begin{center}
$ 0\rightarrow Keru \rightarrow (A\bowtie^{f}J)^{I} \bigoplus (A\bowtie^{f}J)^{p} \rightarrow M\bowtie^{f}J \rightarrow 0 $
\end{center}
where
\begin{eqnarray}
\nonumber u((\alpha_{i},f(\alpha_{i})+j_{i})_{i \in I},(\beta_{i},f(\beta_{i})+k_{i})_{i=1}^{i=p}) &=& \sum_{i \in I}(\alpha_{i},f(\alpha_{i})+j_{i})(m_{i}, f(m_{i})) \\
\nonumber &+& \sum_{i=1}^{i=p}(\beta_{i},f(\beta_{i})+k_{i})(0, g_{i}) \\
\nonumber &=& \sum_{i \in I}(\alpha_{i}m_{i},f(\alpha_{i}m_{i}))+\sum_{i=1}^{i=p}(0,f(\beta_{i}) g_{i})
\end{eqnarray}
$Keru = U\bowtie^{f^{I}}J^{I} \bigoplus (M\bowtie^{f}J)^{p}$. Where $U=\{(\alpha_{i})_{i \in I} \in A^{I} \diagup \sum_{i \in I}\alpha_{i}m_{i}=0 \}$ (since $\{(g_{i})_{i=1}^{i=p}\}$ be a minimal generating set of $J$). Therefore, we have the isomorphism of $(A\bowtie^{f}J)$-modules
\begin{center}
$M\bowtie^{f}J \cong \frac{(A\bowtie^{f}J)^{I}}{U\bowtie^{f^{I}}J^{I}} \bigoplus \frac{(A\bowtie^{f}J)^{p}}{(M\bowtie^{f}J)^{p}}$
\end{center} It
follows that
\begin{center}
$pd_{A\bowtie^{f}J}(\frac{A\bowtie^{f}J}{M\bowtie^{f}J}) \leq pd_{A\bowtie^{f}J}(M\bowtie^{f}J)$ \hspace{2 cm} (2)
\end{center}
Clearly, (1) and (2) force $pd_{A\bowtie^{f}J}(M\bowtie^{f}J)$ to be infinite.
Now the exact sequence of $(A\bowtie^{f}J)$-modules
\begin{center}
$ 0\rightarrow M\bowtie^{f}J \rightarrow A\bowtie^{f}J \rightarrow \{0\} \times J \rightarrow 0 $
\end{center}
where $v(\alpha,f(\alpha)+j) = (\alpha,f(\alpha)+j)(0, k) = (0, f(\alpha)k)$, easily yields $pd_{A\bowtie^{f}J}(\{0\} \times J) = \infty$, as desired. \\
\qed

\bigskip
\noindent {\bf Proof of Theorem \ref{thm3}} \par
{\bf (1)} Let $K_{3}$ be a $3$-finitely presented $(A\bowtie^{f}J)$-module, and let $\{(k_{i})_{i=1}^{i=p} \}$ be a minimal generating set of $K_{3}$. We want to show that $K_{3}$ is a projective ($A\bowtie^{f}J$)-module. For this, consider the exact sequence of ($A\bowtie^{f}J$)-modules:
\begin{center}
$ 0\rightarrow Keru (=K_{2}) \rightarrow (A\bowtie^{f}J)^{p} \rightarrow  K_{3} \rightarrow 0 $
\end{center}
where $u((\alpha_{i})_{i=1,.....,p}) = \sum_{i=1}^{i=p}\alpha_{i}k_{i}$. Then $K_{2} \subset (M\bowtie^{f}J)^{p}$ since $\{(k_{i})_{i=1}^{i=p} \}$ is a minimal generating set of $K_{3}$, $A\bowtie^{f}J$ is a local ring with maximal ideal $M\bowtie^{f}J$. We prove that $K_{2}=0$. Otherwise, $K_{2} \neq 0$. Let $(m_{i},f^{p}(m_{i})+k_{i})_{\{i=1,...,r\}}$ be a minimal generating set of $K_{2}$, where $m_{i} \in M^{p}, k_{i} \in J^{p}$ for each $i=1,.......,r$. Consider the exact sequence of $(A\bowtie^{f}J)$-modules
\begin{center}
$ 0\rightarrow Kerv (=K_{1}) \rightarrow (A\bowtie^{f}J)^{r} \rightarrow K_{2}  \rightarrow 0 $
\end{center}
where $v((\alpha_{i},f(\alpha_{i})+j_{i})_{i=1,.....,r}) = \sum_{i=1}^{i=r}(\alpha_{i}m_{i},(f(\alpha_{i})+j_{i})(f^{p}(m_{i})+k_{i}))=(\sum_{i=1}^{i=r}\alpha_{i}m_{i},\sum_{i=1}^{i=r}f(\alpha_{i})(f^{p}(m_{i})+k_{i}))$. But $K_{1} \subset (M\bowtie^{f}J)^{r}$. Hence $K_{1}  = U\bowtie^{f^{r}}J^{r}$. Where $U=\{(\alpha_{i})_{i=1,...,r} \in A^{r} \diagup \sum_{i=1}^{i=r}\alpha_{i}m_{i}=0 \}$. Since $K_{3}$ is a $3$-presented $(A\bowtie^{f}J)$-module, then $K_{1}$ is a finitely presented $(A\bowtie^{f}J)$-module. By Lemma \ref{lem1}, $U$ is a finitely generated $A$-module and $J$ is a finitely generated ideal of $f(A)+J$. Let $\{(u_{i})_{i=1}^{i=s} \}$ be a set of generators of $U$, and let $\{(g_{i})_{i=s+1}^{i=s+t} \}$ be a minimal generating set of $J^{r}$, where $u_{i} \in M^{r}, g_{i} \in J^{r}$. Consider the exact sequence of $(A\bowtie^{f}J)$-modules
\begin{center}
$ 0\rightarrow Kerw (=K_{0}) \rightarrow (A\bowtie^{f}J)^{s+t} \rightarrow K_{1}  \rightarrow 0 $
\end{center}
where $w((\alpha_{i},f(\alpha_{i})+j_{i})_{i=1,.....,s+t}) = \sum_{i=1}^{i=s}(\alpha_{i},f(\alpha_{i})+j_{i})(u_{i},f^{r}(u_{i}))+\sum_{i=s+1}^{i=s+t}(\alpha_{i},f(\alpha_{i})+j_{i})(0,g_{i})= \sum_{i=1}^{i=s}(\alpha_{i}u_{i},f^{r}(\alpha_{i}u_{i}))+\sum_{i=s+1}^{i=s+t}(0,f(\alpha_{i})g_{i}))=(\sum_{i=1}^{i=s}
\alpha_{i}u_{i},\sum_{i=1}^{i=s}f^{r}(\alpha_{i}u_{i})+\sum_{i=s+1}^{i=s+t}f(\alpha_{i})g_{i})$. It follows that
$K_{0} \cong W\bowtie^{f^{s}}J^{s} \bigoplus M^{t}\bowtie^{f^{t}}J^{t}$, where
$W=\{(\alpha_{i})_{i=1,...,s} \in A^{s} \diagup \sum_{i=1}^{i=s}\alpha_{i}u_{i}=0 \}$ (since $\{(g_{i})_{i=s+1}^{i=s+t} \}$ is a
minimal generating set of $J^{r}$). By the above sequence $K_{0}$ is a finitely generated $(A\bowtie^{f}J)$-module. So $M$ is a
finitely generated ideal of $A$, the desired contradiction.\\

{\bf (2) (a)} Assume that M contains a regular element $m$ and $J$ is a finitely generated ideal of $(f(A)+ J)$. We must show that $A\bowtie^{f}J$ is not a $(2, d)$-ring, for each integer $d \geqslant 0$. Let $K = (A\bowtie^{f}J)(m,f(m))$ and consider the exact sequence of $(A\bowtie^{f}J)$-modules
\begin{center}
$ 0\rightarrow Kerv \rightarrow A\bowtie^{f}J \rightarrow K \rightarrow 0 $ \hspace{2 cm} $(\ast)$
\end{center}
where $v(\alpha,f(\alpha)+j) = (\alpha,f(\alpha)+j)(m,f(m))=(\alpha m, f(\alpha m)) $. Clearly $Kerv= \{0\} \times J$ that is a finitely generated ideal of $(A\bowtie^{f}J)$ and hence $K$ is a finitely presented ideal of $(A\bowtie^{f}J)$ by a sequence $(\ast)$. On the other hand, $pd_{A\bowtie^{f}J}(Ker(v)) = pd_{A\bowtie^{f}J}(\{0\} \times J) = \infty$ by Lemma \ref{lem3}. So $pd_{A\bowtie^{f}J}(K) = \infty$.
Finally, the exact sequence of $(A\bowtie^{f}J)$-modules
\begin{center}
$ 0\rightarrow K \rightarrow A\bowtie^{f}J \rightarrow \frac{A\bowtie^{f}J}{K} \rightarrow 0 $
\end{center}
yields $\frac{A\bowtie^{f}J}{K}$ a $2$-presented ($A\bowtie^{f}J$)-module with infinite projective dimension
(i.e., $A\bowtie^{f}J$ is not a $(2, d)$-ring, for each $d \geqslant 0$).\\

{\bf (b)} The same reasoning as in the proof of Corollary \ref{cor3} shows that $C$ is a $(3,d)$-ring which is neither
a $(2,d)$-ring $(d\geq 0)$ nor a $(3,d-1)$-ring $(d\geq 1)$, and this completes the proof of Theorem \ref{thm3}.\\
\qed

\bigskip
The following Corollary is an immediate consequence of Theorem
\ref{thm3}. \\

\begin{cor}\label{cor4}
Let $(A, M)$ be a local domain such that $M$ is not finitely generated ideal of $A$, $f: A\rightarrow B$ be a ring homomorphism, and let $J$ be a proper ideal of $B$ such that $J$ is a finitely generated ideal of $f(A)+ J$, $f(M)J=0$ and $J^{2}=0$. Then:
\begin{enumerate}
\item $A\bowtie^{f}J$ is a $(3,0)$-ring which is not a $(2,0)$-ring.
\item Let $A_{1}=A\bowtie^{f}J$, $d$ be an integer, $A_{2}$ be a Noetherian ring of global dimension $d$, and let $C=A_{1} \times A_{2}$ the direct product of $A_{1}$ and $A_{2}$. The ring $C$ is a $(3,d)$-ring which is neither a $(2,d)$-ring nor a $(3,d-1)$-ring.
\end{enumerate}
\end{cor}
\bigskip
Theorem \ref{thm3} enriches the literature with new examples of
$(3,d)$-rings which are neither a $(2,d)$-ring $(d\geq 0 )$ nor a
$(3,d-1)$-ring $(d \geq 1)$, as shown below. \\

\begin{exam}
Let $k$ be any field and $X_{1},X_{2},...,X_{n},... $ be
indeterminate over $K$. Let $A=K[[X_{1},X_{2},..,X_{n},..]]$  the
power series ring in infinite variables over $K$, and Let $M$ be
its maximal ideal. Set $B:= A/M^{2}$, $J:=Bf(m)$ be an ideal of
$B$, where $m \in M$ such that $f(m) \neq 0$, and consider the
canonical ring homomorphism $f: A\rightarrow B$ ($f(x) =
\overline{x}$). Let $C=\mathbb{Z}[X_{1},X_{2},...,X_{d-1}]$, where
$d \geq 1$. Then by Corollary \ref{cor4}, $(A\bowtie^{f}J) \times
C$ is a $(3,d)$-ring which is neither a $(2,d)$-ring nor a
$(3,d-1)$-ring.
\end{exam}
\bigskip
\begin{exam}
Let $k$ be a field and let $A=K[[X]]=K+M$, where $M=XA$. Set $B:= A/M^{2}$, $J:=Bf(m)$ be an ideal of $B$, where $m \in M$ such that $f(m) \neq 0$, and consider the canonical ring homomorphism $f: A\rightarrow B$ ($f(x) = \overline{x}$). Then $A\bowtie^{f}J$ is not an $(n,d)$-ring, for any integers $n, d \geq 0$.
\end{exam}
\proof $A\bowtie^{f}J$ is Noetherian ring since $A$ is by \cite[Proposition 5.6]{AFF1} and by Lemma \ref{lem3}, $pd_{A\bowtie^{f}J}(\{0\} \times J) = \infty$, whence $gldim(A\bowtie^{f}J)=\infty$. Then by \cite[Theorem 1.3]{Co}, $A\bowtie^{f}J$ is not an $(n,d)$-ring, for any integers $n, d \geq 0$, as desired.\\
\qed

\end{section}

\bigskip



\bigskip
\end{section}

\end{document}